\documentclass{article}
\usepackage{amssymb,amsmath}

\newtheorem{proposition}{Proposition}[section]
\newtheorem{theorem}{Theorem}[section]
\newtheorem{definition}{Definition}[section]
\newtheorem{remark}{Remark}[section]
\newtheorem{corollary}{Corollary}[section]
\begin{document}
\begin{center}
\textbf{AN ELEMENTARY METHOD OF CALCULATING AN EXPLICIT FORM OF YOUNG MEASURES IN SOME SPECIAL CASES}
\end{center}
\begin{center}
\textbf{Piotr Pucha{\l}a}\\
%$^{\rm a}$$^{\ast}$\thanks{$^\ast$
{\small\mbox{Institute of
 Mathematics, Czestochowa University of Technology,
 al. Armii Krajowej 21}, 42-200 Cz\c{e}stochowa, Poland\\
Email: piotr.puchala@im.pcz.pl, p.st.puchala@gmail.com}
\end{center}

%\maketitle

\begin{abstract}
We present an elementary method of explicit calculation
of Young measures for certain class of functions. This class contains
in particular functions of a highly oscillatory nature which appear in
optimization problems and homogenization theory. In engineering such 
situation occurs for instance in nonlinear elasticity (solid-solid
phase transition in certain elastic crystals). Young measures
associated with oscillating minimizing sequences gather information about their
oscillatory nature and therefore about underlying
microstructure. The method presented in the paper makes no use of
functional analytic tools. There is no need to use gene\-ralized 
version of the Riemann -- Lebesgue lemma and to calculate weak$^{\ast}$
limits of functions. The main tool is the change of variable
theorem. The method applies both to sequences of periodic and
nonperiodic functions.\\

\bigskip
\end{abstract}

\section[]{Introduction}
\quad Young measures were first introduced by
Laurence Chisholm Young in the paper \cite{Young} in 1937 where he considered
variational problems that do not admit classical solutions. In these
cases minimizing sequences have a highly oscillatory nature and do not converge strongly, but usually $\text{weakly}^{\ast}$ to some function,
which is not the solution of the problem. The idea
of Young was to enlarge function spaces to the measure spaces and to
consider objects which he called ,,generalized trajectories''. In
these measure spaces the minimizing sequences have in some sense
generalized limits -- the Young measures. They are suitable tool to
analyze the oscillatory properties of the minimizing sequences and
further -- the microstructure arising for example from phase
transitions in certain elastic crystals.

We refer to \cite{Attouch, Balder, Carstensen}, 
\cite{Malek, Muller, Pedregal}, \cite{Roubicek}  and the original
papers cited there for further information about Young measures and
their applications in optimization theory, nonlinear elasticity,
numerical analysis, economics and other areas of engineering and
mathematics. 

In this paper we propose a separation-of-variables-like
method of direct calculation of Young measures. First, in Section 2,
we present an introduction justifying our approach to Young
measures. The readers who are not familiar with Young
measures may read this part as a short
introduction to the theory. Section 3 begins with the definition of
quasi-Young measures. Then we formulate propositions and theorem which
enable us to calculate them explicitly avoiding not handy to deal with tools
from functional analysis. In fact, we use only the change of variable
formula. Further we use this theorem to calculate quasi-Young measures
associated with the sequences which appear in various problems in optimization
and homoge\-nization theory. This is the content of Sections 3-5 and the
reader interested in applications only can concentrate on them omitting
Section 2. In the last part we show, that in many cases the
quasi-Young measures associated with the functions and sequences of
functions are equal to the
Young measures associated with them. This means that 
calculation of Young measures describing oscillations of the
minimizing sequences appearing in applications can frequently be significantly
simplified.

\section[]{Theoretical foundations} This part is entirely based on
\cite{Roubicek}, where general formulation and detailed proofs of all the
results of this section can be found.

Let $\varOmega\subset{\mathbb{R}^n}$ be an open, bounded set and let
$K\subset{\mathbb{R}^m}$ be a compact set. Let $U$ denote the set of all
Lebesgue measurable bounded functions $u\colon\varOmega\rightarrow K$.
Then $U$ is the subset of the space $L^{\infty}(\varOmega,{\mathbb{R}^m})$. By
$\textnormal{Car}(\varOmega ,K)$ we will denote the space of
Carath\'eodory functions. Recall that the function
$h\colon\varOmega\times K\rightarrow{\mathbb{R}}$ is called a
Carath\'eodory function if $h(\cdot ,k)$ is measurable for all 
$k\in K$  and $h(x,\cdot )$ is continuous for almost all (a.a.) $x$
(with respect to the Lebesgue measure).We will say that $h_1=h_2$,
$h_1,h_2\in \textnormal{Car}(\varOmega ,K)$, iff for
a.a. $x\in\varOmega$ we have 
$\Vert h_1(x,\cdot )-h_2(x,\cdot )\Vert_{C(K)}=0$, where
$\Vert\cdot\Vert_{C(K)}$ denotes the usual supremum norm in the space
of continuous functions. We endow the space 
$\textnormal{Car}(\varOmega ,K)$ with the norm
$\Vert h\Vert_{\textnormal{Car},K}:=
\int\limits_{\varOmega}\sup\limits_{k\in K}\vert h(x,k)\vert dx.$ 
We can identify this space with the Banach space of Bochner integrable
functions from $\varOmega$ to $C(K)$. Namely, we have
\begin{proposition}\label{isomorph_carath}
The mapping $h\mapsto\hat{h}$ defined by $\hat{h}(x):=h(x,\cdot)$
is the isometric isomorphism between the
spaces $\textnormal{Car}(\varOmega ,K)$ and $L^1(\varOmega ,C(K))$.
\end{proposition}
Consider now imbedding 
$i\colon U\rightarrow L^1(\varOmega ,C(K))^\ast$ defined by
\[
\langle i(u),h\rangle :=\int\limits_{\varOmega}h(x,u(x))dx.
\]
The symbol $Y(\varOmega ,K)$ will stand for the weak$^\ast$ closure of
the set $i(U)$ in $L^1(\varOmega ,C(K))^\ast$. So we have
\[
Y(\varOmega ,K):=\Bigl\{L^1(\varOmega ,C(K))^\ast\ni\eta :\exists 
(u_n)\subset U:i(u_n)\xrightarrow[n\rightarrow\infty]{w^\ast}\eta\;
\textnormal{in}\; L^1(\varOmega ,C(K))^\ast\Bigr\}.
\]
\begin{theorem}
The triple
$\bigl( Y(\varOmega ,K),L^1(\varOmega ,C(K))^\ast ,i\bigr)$ is the
convex com\-pac\-ti\-fi\-cation of the set $U$.
\end{theorem}
Recall that by the Riesz theorem we have
$C(K)^\ast =\textnormal{rca}(K)$, where $\textnormal{rca}(K)$ is the
set of regular (signed) measures on $K$. We say that the mapping
$\nu\colon\varOmega\rightarrow\textnormal{rca}(K)$ is weakly
measurable if for all $z\in C(K)$ the mapping 
$x\mapsto\langle\nu (x),z\rangle$ is measurable. We often write
$\nu_x$ instead of $\nu (x)$. The set of all such mappings $\nu$
fulfiling additionally the condition
\[
\Vert\nu\Vert_{L_{w}^{\infty}(\varOmega ,\textnormal{rca}(K))}:=
\textnormal{ess}\sup\bigl\{\Vert\nu (x)\Vert_{\textnormal{rca}(K)}:
x\in\varOmega\bigr\}<+\infty
\]
is a normed linear space $L_{w}^{\infty}(\varOmega ,\textnormal{rca}(K))$. 
By the Dunford -- Pettis theorem this space is isometrically
isomorphic to the the space $L^1(\varOmega ,C(K))^\ast$ (see Section
1.4 of \cite{Roubicek}). Namely, we
have the following proposition.
\begin{proposition}\label{young_measures}
Let $h\in L^1(\varOmega ,C(K))$. Define the mapping 
$\eta\in L^1(\varOmega ,C(K))^\ast$ by
\[
\langle\eta ,h\rangle :=\int\limits_{\varOmega}
\Bigl(\int\limits_{K}
h(x,k)d\nu_x(k)\Bigr)dx.
\]
Then the mapping
\[
\psi\colon L^{\infty}_w(\varOmega,\textnormal{rca}(K))\ni
\nu\mapsto\psi (\nu ):=\eta\in L^1(\varOmega ,C(K))^\ast
\]
realizes the isometric isomorphism between the spaces
$L^{\infty}_w(\varOmega,\textnormal{rca}(K))$
and $L^1(\varOmega ,C(K))^\ast$.
\end{proposition}
In what follows the symbol $\textnormal{rca}^1(K)$ will denote the
subset of $\textnormal{rca}(K)$ consisting of all probability measures
on $K$. In particular the Dirac measure $\delta_k$ concentrated in
$k\in K$ belongs to $\textnormal{rca}^1(K)$.

Define
\[
\mathcal{Y}(\varOmega ,K):=\bigl\{\nu=(\nu (x))\in
L^{\infty}_w(\varOmega,\textnormal{rca}(K)):\nu (x)\in
\textnormal{rca}^1(K)\;\textnormal{for a.a }x\in\varOmega\bigr\}.
\]
This is the set of all Young measures. Adopting the usual notation we will
write $\nu_x$ or $(\nu_x)_{x\in\varOmega}$ instead of $\nu (x)$.

Define now $\delta\colon U\ni u\mapsto \delta_{u(x)}\in
\mathcal{Y}(\varOmega ,K)$.
\begin{theorem}
The function $\psi$ defined in Proposition \ref{young_measures} maps
the set $\mathcal{Y}(\varOmega ,K)$ onto the set $Y(\varOmega ,K)$.
Moreover, the convex compactifications
$\bigl( Y(\varOmega ,K), L^1(\varOmega ,C(K))^\ast ,i\bigr)$ 
and $\bigl(\mathcal{Y}(\varOmega ,K),
L^{\infty}_w(\varOmega,\textnormal{rca}(K)),\delta\bigr)$ of the set
$U\subset L^{\infty}(\varOmega,{\mathbb{R}^m})$ are equivalent. 
\end{theorem}
This means that for every function $u\in U$ there exists an element of
$Y(\varOmega ,K)$ -- this is the Young measure associated with $u$.

\section[]{Quasi-Young measures}
In this part $\varOmega\subset\mathbb{R}^n$ is an open, bounded and
Lebesgue measurable set of measure $M$, $d\mu(x) :=\tfrac{1}{M}dx$ with $n$ -- 
dimensional Lebesgue measure $dx$. Further, $K\subset\mathbb{R}^n$ is a compact set,
$u\colon\varOmega \rightarrow K$ is a measurable function such that  
$\overline{u(\varOmega )}=K$ and $\beta\in C(K,\mathbb{R})$. If $u$
has partial derivatives for $x\in\varOmega$,
determinant of the Jacobi matrix of $u$ is denoted by $J_u$, that is
\[
J_u(x)=\det\Bigl[\tfrac{\partial u_i}{\partial x_j}(x)\Bigr].
\]
\begin{definition}\label{quasi_def}
We say that a family of probability measures 
$\nu =(\nu_x)_{x\in\varOmega}$ is a quasi-Young measure associated with
the measurable function $u\colon\varOmega\rightarrow K$, if for every
continuous function  
$\beta\colon K\rightarrow\mathbb{R}$ there holds an equality
\begin{equation}\label{quasi_def_eq}
\int\limits_{K}\beta(k)d\nu_x(k)=
\int\limits_{\varOmega}\beta(u(x))d\mu(x).
\end{equation}
\end{definition}
\begin{proposition}\label{mon_quasi}
Let $u$ be a function such that its inverse $u^{-1}$ is continuously
differentiable. Then a quasi-Young measure associated with $u$ is a
measure 
that is absolutely continuous with respect to the Lebesgue measure on
$K$. Its density is equal to 
$\tfrac{1}{M}\vert J_{u^{-1}}\vert$.
\end{proposition}
{\itshape Proof.}
Using the change of variable theorem we get
\[
\int\limits_K\beta (k)d\nu_x(k)=\int\limits_{\varOmega}
\beta (u(x))d\mu(x)=
\int\limits_K\beta (y)\tfrac{1}{M}\vert J_{u^{-1}}(y)\vert dy.
\]
\begin{flushright}
$\square$
\end{flushright}
\begin{remark}
We see that the quasi-Young measure does not depend on the va\-riab\-le
$x$. This is a homogeneus quasi-Young measure.
\end{remark}
Recall that for any set $A$ the symbol $\chi_{A}$ denotes the
characteristic function of $A$, i.e.
\[
\chi_{A}(x)=
\begin{cases}
0, & x\notin A\\
1, & x\in A.
\end{cases}
\]
We will now introduce a partition of $\varOmega$ 
into $n$ open subsets $\varOmega_1,\dots, \varOmega_n$ such that
\begin{enumerate}
\item[(a)]
$\varOmega_i\cap \varOmega_j=\emptyset\textnormal{ for }i\neq j$;
\item[(b)]
$\bigcup\limits_{i=1}^n\overline{\varOmega}_i=\overline{\varOmega}$;
\end{enumerate}
In the sequel the symbol $\{\varOmega\}$ will denote the partition
of $\varOmega$ defined above.
\begin{proposition}\label{gestosc_klejona}
Let the function $u\colon\varOmega\rightarrow K$ has the form
\[
u:=\sum\limits_{i=1}^{n}u_i\chi_{\varOmega_i},
\]
where $u_i$ satisfies the assumptions of the Proposition
\ref{mon_quasi}, $\overline{u_i(\varOmega )}=K$
and $\varOmega_i$ belongs to $\{\varOmega\}$, $i=1,2,\dots,n$. Then the
quasi-Young measure associated with $u$ is absolutely continuous with
respect to the Lebesgue measure $dy$ on $K$ with the density
\begin{equation}\label{gestosc_klejona_eq}
g=\tfrac{1}{M}\sum\limits_{i=1}^{n}\vert J_{u_i^{-1}}(y)\vert .
\end{equation}
\end{proposition}
{\itshape Proof} We proceed by induction with respect to the number of
,,components'' $u_i$ of $u$. The case $n=1$ is Proposition 
\ref{mon_quasi}. Suppose that the statement is true for some integer
$l>1$. Then 
\[ 
g_l=\tfrac{1}{M}\sum_{i=1}^{l}\vert J_{u_i^{-1}}(y)\vert ,
\]
so we have
\begin{align}
\int\limits_{K}\beta (k)d\nu_x(k)=\int\limits_{\varOmega}
\beta (u(x))d\mu(x)=\int\limits_{\varOmega}\beta\biggl(
\sum\limits_{i=1}^{l+1}u_i(x)\chi_{\varOmega_i}(x)\biggr)d\mu(x) =\notag \\
=\sum\limits_{i=1}^{l+1}\int\limits_{\varOmega_i}\beta (u_i(x))d\mu(x) =
\int\limits_K\beta (y)\tfrac{1}{M}
\sum\limits_{i=1}^{l+1}\vert J_{u_i^{-1}}(y)\vert dy.\notag
\end{align}
\begin{flushright}
$\square$
\end{flushright}
We can prove analogous result when function $u$ is 'built' of
functions having continuously differentiable inverses.
\begin{proposition}\label{gestosc_klejona_monot}
Let the continuous function $u\colon\varOmega\rightarrow K$ has the form
\[
u:=\sum\limits_{i=1}^{n}u_i\chi_{\varOmega_i},
\]
where for each $i=1,2,\dots, n$ $\varOmega_i\in\{\varOmega\}$ and 
the function $u_i$ has continuously differentiable inverse $u_i^{-1}$.
Assume further that the set $\overline{u_i(\varOmega_i )}$ is compact,
$i=1,\dots, n$,
$\bigcup\limits_{i=1}^n\overline{u_i(\varOmega_i)}=K$ and that 
$u_i(\varOmega_i)\cap u_j(\varOmega_j)=\emptyset$ for $i\neq j$.
Then the quasi-Young measure associated with $u$ is absolutely continuous with
respect to the Lebesgue measure $dy$ on $K$ with the density
\begin{equation}\label{gestosc_klejona_monot_eq}
g=\tfrac{1}{M}
\sum\limits_{i=1}^{n}\vert J_{u_i^{-1}}(y)\vert\cdot\chi_{u_i(\varOmega_i)} .
\end{equation} 
\end{proposition}
Analogous results hold for $u$ constant or piecewise constant.
\begin{proposition}\label{const_quasi}
\begin{enumerate}
\item[(a)]
assume that $u$ is a constant function: $\forall x\in\varOmega$ $u(x)=p$,
$p$ -- a fixed vector in $\mathbb{R}^n$. Then the quasi-Young measure
associated with $u$ is the Dirac measure $\delta_p$.
\item[(b)]
let the function $u$ has the form
\[
u:=\sum\limits_{i=1}^{n}p_i\chi_{\varOmega_i},
\]
where $\varOmega_i\in\{\varOmega\}$ and 
$p_i$ are fixed vectors in $\mathbb{R}^n$, $i=1,\dots ,n$. Then its quasi-Young
measure is a convex combination of Dirac measures:
\[
\nu_x=\tfrac{1}{M}
\sum\limits_{i=1}^{n}m_i\delta_{p_i},
\]
where $m_i$ is the Lebesgue measure of $\varOmega_i\in\{\varOmega\}$,
$i=1,2,\dots ,n$.
\end{enumerate}
\end{proposition}
{\itshape Proof.} It is enough to prove (a). We have
\[
\int\limits_{K}\beta (k)d\nu_x(k)=\int\limits_\varOmega\beta(u(x))d\mu(x)=
\beta (p)=\int\limits_{K}\beta (y)d\delta_p.
\]
\begin{flushright}
$\square$
\end{flushright}

Consider now a family $\Pi^{(k)}$ of open partitions of $\varOmega$,
such that $\Pi^{(1)}=\varOmega$, 
$\Pi^{(k)}=\bigl\{\varOmega_i^{(k)}\bigr\}_{i=1}^{k}$, $k=2,3,\dots$,
and for any fixed $2\leq k\in\mathbb{N}$ we have:
\begin{itemize}
\item[(a)]
$\varOmega_i^{(k)}\cap\varOmega_j^{(k)}=\emptyset$ for $i\neq j$, 
$1\leq i,\, j\leq k;$
\item[(b)]
$\bigcup\limits_{i=1}^{k}\overline{\varOmega_i^{(k)}}=
\overline{\varOmega};$
\item[(c)]
$\forall i\in\{1,\dots,k\},\;\varOmega_i^{(k)}=
\textnormal{interior}\overline{(\varOmega_i^{(k)}}).$
\end{itemize} 
Consider further a family $\{u^{(k)}\}$ of Lebesgue measurable bounded
functions from $\varOmega$ to $K$, associated with $\Pi^{(k)}$, 
$k=1,2,\dots$, in such a way that:
\begin{itemize}
\item[(i)]
$u^{(1)}=u$;
\item[(ii)]
for fixed $2\leq k\in\mathbb{N}$ we have
\begin{equation}
u^{(k)}(x)=
\begin{cases}
u_1^{(k)}(x), & x\in\varOmega_1^{(k)}\\
u_2^{(k)}(x), & x\in\varOmega_2^{(k)}\\
\cdot \\
\cdot \\
\cdot \\
u_k^{(k)}(x), & x\in\varOmega_k^{(k)}.\\
\end{cases}
\end{equation}
\end{itemize}
In view of Propositions \ref{gestosc_klejona} and
\ref{const_quasi} the following theorem is true.
\newpage
\begin{theorem}\label{generating}
Consider the sequence $(u^{(k)})$ of functions associated with the
family $\Pi^{(k)}$ of partitions of $\varOmega$ such that for any
fixed $k\in\mathbb{N}$ we have
\begin{itemize}
\item[(a)]
for $i=1,\dots ,k$ the functions $u_i^{(k)}\colon\varOmega\to K$ 
have continuously differentiable inverses $[u_i^{(k)}]^{-1}$,
$\overline{u_i^{(k)}(\varOmega )}=K$ and moreover
\[
\tfrac{1}{M}\sum\limits_{i=1}^{k}\vert J_{[u_i^{(k)}]^{-1}}(y)\vert
=\textnormal{const}=:g;
\]
\end{itemize}
or
\begin{itemize}
\item[(b)]
for $i=1,\dots ,k$ the functions  $u_i^{(k)}$ are constant on
$\varOmega_i^{(k)}$ with value $p_i^{(k)}$.
\end{itemize}
Then the quasi-Young measure associated with each term $u^{(k)}$ is
respectively: 
\begin{itemize}
\item[(a)]
absolutely continuous with respect to the Lebesgue measure $dy$ on $K$
with the density $g$;
\end{itemize}
or
\begin{itemize}
\item[(b)]
convex combination of Dirac measures:
\[
\nu_x=\tfrac{1}{M}
\sum\limits_{i=1}^{k}m_i^{(k)}\delta_{p_i^{(k)}},
\]
where $m_i^{(k)}$ is the Lebesgue measure of
$\varOmega_i^{(k)}\in\Pi^{(k)}$.
\end{itemize}
\end{theorem}
\section[]{Sequences of oscillating functions}
\quad Rapidly oscillating functions play an important role in 
homogenization theory and in the optimization problems when we are
seeking minima of integral functionals with nonconvex integrands. In
this second case minimizing sequences do not converge strongly in an
appropriate topology, but usually the convergence is only the weak$^\ast$
one. Finer and finer oscillations of the elements of the minimizing
sequence around its weak$^\ast$ limit give some information about
microstructure. The Young measures associated with the minimizing
sequences 'capture' these oscillations.

\mbox{Let $(c_n)$ be a monotonically increasing sequence of natural numbers:
$\lim\limits_{n\rightarrow\infty}c_n=+\infty$.}
\mbox{Consider Lebesgue measurable function 
$u\colon\varOmega :=]0,1[\to K:=[0,1]$ with 
$\overline{u(\varOmega)}=K$} and
the sequence $(u_n)$ of functions defined by
$u_n(x):=u(c_nx)$, $n\in\mathbb{N}$. We will call 
the function $u$ ''the function ge\-ne\-rating fast oscillating
sequence''. For each $n\in\mathbb{N}$ we can write
\[
u_n(x)=
\begin{cases}
u(c_nx), & x\in ]0,\tfrac{1}{c_n}]\\
u(c_nx-1), & x\in ]\tfrac{1}{c_n},\tfrac{2}{c_n}]\\
\cdot \\
\cdot \\
\cdot \\
u(c_nx-(c_n-1)), &x\in ]\tfrac{c_n-1}{c_n}, 1[.
\end{cases}
\]
The following proposition is in fact a special case of 
Theorem \ref{generating}.
\begin{proposition}\label{quasi_osc_const}
\begin{enumerate}
\item[(a)] 
let $u$ be a function strictly monotonic and differentiable on $\varOmega$.
Then the quasi-Young measure associated with the function $u_n$,
$n\in\mathbb{N}$, is equal to the quasi-Young measure associated with
the function $u$ generating the sequence $(u_n)$;
\item[(b)] 
an analogous result holds if the generating function $u$ is
piecewise constant on $\varOmega$.
\end{enumerate}
\end{proposition}
{\itshape Proof}
It is enough to prove (a). Let $u$ be a strictly monotonic
differentiable function. We can assume 
that $u$ is increasing. Choose and fix $n\in\mathbb{N}$. Then we have
\begin{align}
\int\limits_{0}^{1}\beta(k)d\nu_x(k)=\int\limits_{0}^{1}\beta(u_n(x))dx=
\sum\limits_{i=1}^{c_n}\int\limits_{\tfrac{i-1}{c_n}}^{\tfrac{i}{c_n}}
\beta(u(c_nx-(i-1))dx= \notag \\
=\sum\limits_{i=1}^{c_n}\frac{1}{c_n}\int\limits_{0}^{1}\beta(u(t))dt=
\int\limits_{0}^{1}\beta(y)(u^{-1})'(y)dy. \notag
\end{align}
An inductive argument completes the proof.
\begin{flushright}
$\square$
\end{flushright}
%\newpage
\begin{corollary}\label{quasi_osc_seq} 
\begin{enumerate}
\item[(a)] 
let the function $u\colon\varOmega\rightarrow K$ generating the fast
oscillating sequence has the form
\[
u:=\sum\limits_{i=1}^{n}u_i\chi_{\varOmega_i},
\]
where $u_i$ is strictly monotonic and differentiable on $\varOmega_i$,
$\overline{u_i(\varOmega_i )}=K$ 
and $\varOmega_i$ belongs to $\{\varOmega\}$, $i=1,2,\dots,n$.
Then for any $n\in\mathbb{N}$ the quasi-Young measure
associated with $u_n$ is the same as the quasi-Young measure associated
with $u$;
\item[(b)] 
an analogous result is true if the generating function $u$ has the form
\[
u:=\sum\limits_{i=1}^{n}p_i\chi_{\varOmega_i},
\]
where $\varOmega_i\in\{\varOmega\}$ and 
$p_i$ are real constants, $i=1,\dots ,n$.
\end{enumerate}
\end{corollary}
\begin{remark}
Let $u$ be a function generating the fast oscillating sequence $(u_n)$.
Observe that  for each $n\in\mathbb{N}$
\[
\int\limits_{0}^{1}u_n(x)dx=\sum\limits_{i=1}^{c_n}
\int\limits_{\tfrac{i-1}{c_n}}^{\tfrac{i}{c_n}}
u(c_nx-(i-1))dx=\int\limits_{0}^{1}u(x)dx,
\]
so geometrically the area between the graph of $u$ and the $x$-axis is the
same as the area between the graph of $u_n$ and the $x$-axis,
$n\in\mathbb{N}$. In this 
case the sequence of quasi-Young measures is constant and hence
trivially $\textnormal{weak}^\ast$ convergent. 
\end{remark}
\section[]{Examples}
\indent\indent We now give some examples of direct computation of
quasi-Young measures. The examples of oscillating sequences are
taken from the items listed in the references (see [\cite{Cioranescu}, \cite{Pedregal} -- \cite{Roubicek}). We can see that the
method described in the above sections gives the same results without
using functional analytic tools. In particular we do not need a
genera\-li\-zed version of the Riemann -- Lebesgue lemma.\\
\\
(a) Let $a,b>0$, $\varOmega =]0,a[$. Let the  function $u$
\[
u(x):=
\begin{cases}
\tfrac{2b}{a}x, & x\in \bigl]0,\tfrac{a}{2}\bigr]\\
-\tfrac{2b}{a}x+2b, & 
x\in\bigl]\tfrac{a}{2},a\bigr[,
\end{cases}
\] 
generates the sequence
\[
u_n(x):=
\begin{cases}
\tfrac{2nb}{a}x-2bk, & x\in \bigl]\tfrac{ak}{n},\tfrac{(2k+1)a}{2n}\bigr]\\
-\tfrac{2nb}{a}x+2b(k+1), & 
x\in\bigl]\tfrac{(2k+1)a}{2n},\tfrac{(k+1)a}{n}\bigr[,
\end{cases}
\]
$n\in\mathbb{N}$ and $k=0,1,\dots ,n-1$. We have
\begin{align}
\int\limits_{0}^{b}\beta (k)d\nu_x(k)=\int\limits_{0}^{a}
\beta(u_n(x))d\mu(x)
=\underbrace{\sum\limits_{k=0}^{n-1}
\int\limits_{\tfrac{ak}{n}}^{\tfrac{(2k+1)a}{2n}}
\beta\bigl(\tfrac{2nb}{a}x-2bk\bigr)d\mu(x)}_\mathrm{I_1}+\notag \\
+\underbrace{\sum\limits_{k=0}^{n-1}
\int\limits_{\tfrac{(2k+1)a}{2n}}^{\tfrac{(k+1)a}{n}}
\beta\bigl(-\tfrac{2nb}{a}x+2b(k+1)\bigr)d\mu(x)}_\mathrm{I_2},\notag
\end{align}
\mbox{where $d\mu(x)=\tfrac{1}{a}dx$. The k-th integral in $I_1$ is equal to 
$\tfrac{1}{2an}\int\limits_{0}^{b}\beta (y)\tfrac{a}{b}dy$, 
$k=0,\dots ,n-1$,} (the same in $I_2$) so finally
\[
\int\limits_{0}^{b}\beta (k)d\nu_x(k)=
\int\limits_{0}^{b}\beta (y)\tfrac{1}{b}dy,
\]
i.e. the quasi-Young measure associated with the function $u_n$ has the
form $\nu_x=\tfrac{1}{b}dy$, $n\in\mathbb{N}$. Thus we have the
constant sequence of 
measures absolutely continuous with respect to the Lebesgue measure
$dy$ on $K=[0,b]$ with the density equal to $\tfrac{1}{b}$. Observe
that for each $n\in\mathbb{N}$
\[
\int\limits_{0}^{a}u_1(x)d\mu(x)=\int\limits_{0}^{a}u_n(x)d\mu(x)
\]
and that calculating the density is in fact reduced to compute the
sum of absolute values 
of the inverses of the slopes of $u_1$ and
dividing this sum by the Lebesgue measure of the set $\varOmega$.\\
(b) Let $\varOmega =]0,1[$ and define the function
\[
u_1(x):=
\begin{cases}
3x, & x\in\bigl]0,\tfrac{1}{6}\bigr]\\
\tfrac{3}{2}x+\tfrac{1}{4}, & x\in\bigl]\tfrac{1}{6},\tfrac{1}{2}\bigr]\\
-\tfrac{3}{2}x+\tfrac{7}{4}, & x\in\bigl]\tfrac{1}{2},\tfrac{5}{6}\bigr]\\
-3x+3, & x\in\bigl]\tfrac{5}{6},1\bigr[.
\end{cases}
\]
Let $(u_n)$ be a sequence of periodic functions with n-th element
built with $n$ such shaped 'teeth' in $\varOmega$. As above, we have
for $n\in\mathbb{N}$ 
$\int\limits_{0}^{1}u_1(x)dx=\int\limits_{0}^{1}u_n(x)dx$, so it is
enough to consider $u_1$. We have
\[
\int\limits_0^1\beta (k)d\nu_x(k)=
\int\limits_0^1\beta (u_1(x))dx=
\int\limits_{0}^{\tfrac{1}{2}}\beta (y)\tfrac{2}{3}dy+
\int\limits_{\tfrac{1}{2}}^{1}\beta (y)\tfrac{4}{3}dy.
\]
and thus $\nu_x =f(y)dy$, where
\[
f(y)=
\begin{cases}
\tfrac{2}{3}, & y\in [0,\tfrac{1}{2}]\\
\tfrac{4}{3}, & y\in ]\tfrac{1}{2},1].
\end{cases}
\]
(c) Let $\varOmega =]0,1[$ and $u_n(x)=\sin (2\pi nx)$. Dividing
$u_1$ into monotonic parts and applying our procedure we get
\begin{multline}
\int\limits_{-1}^{1}\beta (k)d\nu_x(k)=\\
=\int\limits_0^1\beta (y)\frac{dy}{2\pi\sqrt{1-y^2}}+
\int\limits_{-1}^{1}\beta (y)\frac{dy}{2\pi\sqrt{1-y^2}}+
\int\limits_{-1}^{0}\beta (y)\frac{dy}{2\pi\sqrt{1-y^2}}=
\int\limits_{-1}^{1}\beta (y)\frac{dy}{\pi\sqrt{1-y^2}},\notag
\end{multline}
obtaining finally
\[
\nu_x=\frac{dy}{\pi\sqrt{1-y^2}}.
\]
(d) Let $a,b$ be fixed positive real numbers and let the function 
\[
u=
\begin{cases}
a, & x\in ]0,\tfrac{2}{3}]\\
b, & x\in [\tfrac{2}{3},2[
\end{cases}
\]
generate the sequence $(u_n)$. Applying our procedure we get
\[
\int\limits_{K}\beta (k)d\nu_x(k)=
\int\limits_{0}^{\tfrac{2}{3}}\beta (a)d\mu(x)+
\int\limits_{\tfrac{2}{3}}^{2}\beta (b)d\mu(x)=
\int\limits_K\beta (y)(\tfrac{1}{3}\delta_a+\tfrac{2}{3}\delta_b)(dy),
\]
which means that $\nu_x$ is a purely singular measure
\[
\nu_x=\tfrac{1}{3}\delta_a+\tfrac{2}{3}\delta_b.
\]
In this case the coefficients multiplying the Dirac deltas are equal
to the lenght of the respective interval divided by the Lebesgue
measure of $\varOmega$.\\
(e) Let
\[
u_n(x):=
\begin{cases}
(x(n+k-1)-k+1)\tfrac{n+k}{n}, & x\in
]\tfrac{k-1}{n+k-1},\tfrac{k}{n+k}[,\; k\in\mathbb{N}\;\textnormal{odd}\\
(k-x(n+k))\tfrac{n+k-1}{n}, & 
x\in [\tfrac{k-1}{n+k-1},\tfrac{k}{n+k}[,\; k\in\mathbb{N}\;\textnormal{even}.
\end{cases}
\]
This is an example of an oscillating sequence of nonperiodic
functions. Observe that for each $n\in\mathbb{N}$ the sum of the
inverses of absolute values of the slopes of components of $u_n$ is
constant:
\[
n\sum\limits_{k=1}^{\infty}\frac{1}{(k+(n-1))(k+n)}=1,
\] 
which yields the result $\nu_x=1\cdot dx$.\\
(f) Let $\varOmega =]0,1[\times ]0,1[\ni (x_1,x_2)$
and consider the 
sequence of matrix valued functions of the form 
\[
u_n(x):=\chi_{]0,3/4[}(n(x_1+x_2))(1,1)\otimes (1,1).
\]
Denoting by $\langle a\rangle$ the integer part of $a$ we
can write $u_n$ in matrix form
\[
u_n(x)=
\begin{cases}
\left(\begin{array}{cc}
1 & 1\\
1 & 1
\end{array}\right), & 0<n(x_1+x_2)-\langle n(x_1+x_2)\rangle
<\tfrac{3}{4},\\
\\
\left(\begin{array}{cc}
0 & 0\\
0 & 0
\end{array}\right), & \tfrac{3}{4}<n(x_1+x_2)-\langle n(x_1+x_2)\rangle
<1.
\end{cases}
\]
Observe that for each $n\in\mathbb{N}$ the Lebesgue measure
of the set\\
$\varOmega_1:=\{\varOmega\ni (x_1,x_2):
0<n(x_1+x_2)-\langle n(x_1+x_2)\rangle <3/4\}$ equals $3/4$, while the
\mbox{Lebesgue measure of the set
$\varOmega_2:=\{\varOmega\ni (x_1,x_2):
3/4<n(x_1+x_2)-\langle n(x_1+x_2)\rangle <1\}$} equals $1/4$. Denote
by $A$ the matrix with all elements equal to 1 and by $O$ the zero
matrix. The quasi-Young measure associated with $u_n$, $n=1,2,\dots$,
is 
\[
\nu_x=\tfrac{3}{4}\delta_A+\tfrac{1}{4}\delta_O .
\]

\section[]{Quasi-Young measures are Young measures} 
\indent\indent In this section it will be shown that quasi -- Young
measures in the above sections are in fact Young measures. Recall that
$\beta\in C(K,\mathbb{R})$ and that $\{\varOmega\}$ is an open
partition of $\varOmega$  into $n$ open subsets 
$\varOmega_1,\dots,\varOmega_n$ such that 
\begin{enumerate}
\item[(a)]
$\varOmega_i\cap \varOmega_j=\emptyset\textnormal{ for }i\neq j$;
\item[(b)]
$\bigcup\limits_{i=1}^n\overline{\varOmega}_i=\overline{\varOmega}$ .
\end{enumerate}
We will denote by the letter $\alpha$ the continuous function  
$\alpha\colon\overline{\varOmega}\rightarrow\mathbb{R}$. Further, 
$\alpha\otimes\beta$ will stand for the tensor product of functions
$\alpha$ and $\beta$, namely\\ 
$(\alpha\otimes\beta)(x,k):=\alpha (x)\cdot\beta (k)$. Recall that by
Proposition \ref{isomorph_carath} the space 
$\textnormal{Car}(\varOmega ,K)$ of the Carath\'eodory functions is
isometrically isomorphic with the space $L^1(\varOmega ,C(K))$. We
will need Theorem I.5.25 (3) from \cite{Warga}.
According to this theorem the linear hull
$C(\overline{\varOmega})\otimes C(K)$ of the set
$\{f\otimes g:f\in C(\overline{\varOmega}),\, g\in C(K)\}$ is dense in
the space $L^1(\varOmega ,C(K))$.

We assume that $u\colon\varOmega\rightarrow K$ is bounded, Lebesgue
measurable function invertible on $\varOmega$, with
$\overline{u(\varOmega )}=K$. 
\begin{theorem}\label{result}
\begin{enumerate}
\item[(a)] 
let $u$ be
such that its inverse $u^{-1}$
is continuously differentiable on $K$. 
Then the quasi-Young measure associated with $u$ is
equal to the Young measure associated with $u$;
\item[(b)] 
let the function $u$ has the form
\[
u:=\sum\limits_{i=1}^{n}p_i\chi_{\varOmega_i},
\]
where $\varOmega_i\in\{\varOmega\}$ and 
$p_i$ are fixed vectors in $\mathbb{R}^n$, $i=1,\dots ,n$.
Then the quasi-Young measure associated with $u$ is
equal to the Young measure associated with $u$.
\end{enumerate}
\end{theorem}
{\itshape Proof.} It is enough to prove (a). Let $h$ be a Carath\'eodory function, 
$\alpha_{u,\varepsilon}\in C(\overline{\varOmega}, {\mathbb{R}})$ and
$\beta_{u,\varepsilon}\in C(K,{\mathbb{R}})$.
Thanks to Proposition \ref{young_measures} we can write
\begin{align}\notag
\Big\vert\int\limits_{\varOmega}
\Bigl(\int\limits_{K}h(x,k)d\nu_x(k)\Bigr)dx-
\int\limits_{\varOmega}
\Bigl(\int\limits_{K}h(x,k)\vert J_{u^{-1}}(y)\vert dk\Bigr)dx\Big\vert\leq
\\ \notag
\leq
\Big\vert\int\limits_{\varOmega}
\Bigl(\int\limits_{K}h(x,k)d\nu_x(k)\Bigr) dx-
\int\limits_{\varOmega}\alpha_{u,\varepsilon}(x)
\Bigl(
\int\limits_{K}\beta_{u,\varepsilon}(k)d\nu_x(k)\Bigr)dx\Big\vert +
\\ \notag
+\Big\vert\int\limits_{\varOmega}\alpha_{u,\varepsilon}(x)
\Bigl(\int\limits_{K}\beta_{u,\varepsilon}(k)d\nu_x(k)\Bigr)dx
-\int\limits_{\varOmega}\alpha_{u,\varepsilon}(x)
\Bigl(\int\limits_{K}\beta_{u,\varepsilon}(k)\vert J_{u^{-1}}(y)\vert
dk\Bigr)dx\Big\vert+ 
\\ \notag
+\Big\vert\int\limits_{\varOmega}\int\limits_{K}\Bigl( 
h(x,k)-\alpha_{u,\varepsilon}(x)\cdot\beta_{u,\varepsilon}(k)\Bigr)
\vert J_{u^{-1}}(y)\vert dkdx\Big\vert .
\end{align}
Choose and fix the function $u$ and
$\varepsilon >0$. From the denseness results stated above there exists
function
$h_{u,\varepsilon}=\alpha_{u,\varepsilon}\otimes\beta_{u,\varepsilon}$,
such that the first and the third term on the right hand side are 
smaller than $\varepsilon$. The second term vanishes thanks to the
Proposition \ref{mon_quasi}.  
\begin{flushright}
$\square$
\end{flushright}
\begin{corollary}
The quasi-Young measures of theorems of Sections 3 and 4, associated
with the functions $u$ of these theorems, are equal to the Young
measures associated with them.  
\end{corollary}
Theorem \ref{result} enables us to calculate explicit form of Young
measures associa\-ted with functions or to the sequences of functions
considering only the inner integral 
$\int\limits_{K}\beta (k)d\nu_x(k)$. There is no need of periodical
extending $u$ to use the
ge\-ne\-ra\-lized Riemann -- Lebesgue lemma and to calculate 
$\textnormal{weak}^{\ast}$ limits of sequences of functions.\\

\textbf{Acknowledgements.}
Author would like to thank Professor
Zdzis{\l}aw Naniewicz and Kazi\-mierz Br\c{a}giel for many helpful
discussions.

\end{document}